\documentclass[letterpaper,10pt]{IEEEtran}
\usepackage[sort&compress]{natbib}
\usepackage{amsmath,amssymb}
\usepackage[usenames]{color}
\usepackage{xcolor}
\usepackage{caption,subcaption}
\usepackage{graphicx}
\usepackage{enumerate}
\usepackage{nicefrac}

\newtheorem{theorem}{\bf Theorem}
\newtheorem{proposition}{\bf Proposition}
\newtheorem{definition}{\bf Definition}
\newtheorem{lemma}{\bf Lemma}
\newtheorem{corollary}{\bf Corollary}

\newtheorem{remark}{\bf Remark}

\def\ds{\displaystylse}

\newcommand{\dfb}{\stackrel{\Delta}{=}}

\def\E{\mathbb{E}}

\def\calV{\mathcal{V}}

\DeclareMathOperator*{\argmax}{argmax} 

\def\qed{\hfill$\blacksquare$}

\def\one{\mathbf{1}}

\newcommand{\ups}{\upsilon}
\newcommand{\ov}{\overline}
\newcommand{\mc}{\mathcal}
\newcommand{\0}{\mathbf{0}}

\newcommand{\R}{\mathbb{R}}
\renewcommand{\P}{\mathbb{P}}
\renewcommand{\E}{\mathbb{E}}
\newcommand{\ba}{\begin{array}}
\newcommand{\ea}{\end{array}}
\renewcommand{\ds}{\displaystyle}
\newcommand{\beq}{\begin{equation}}
\newcommand{\eeq}{\end{equation}}
\newcommand{\beqn}{\begin{equation*}}
\newcommand{\eeqn}{\end{equation*}}
\newcommand{\be}{\begin{equation}}
\newcommand{\ee}{\end{equation}}
\DeclareMathOperator*{\argmin}{argmin}

\begin{document}

\title{Wisdom of Crowds\\ Through Myopic Self-Confidence Adaptation

\thanks{Giacomo Como and Fabio Fagnani are with the Department of Mathematical Sciences and Anton Proskurnikov is with the Department of Electronics and Telecommunications, both at Politecnico di Torino, Torino, Italy.
 Email: \{{\tt\small{giacomo.como, fabio.fagnani, anton.proskurnikov}\}@polito.it}.  
}

}

\author{Giacomo Como, Fabio Fagnani, Anton Proskurnikov} 

\maketitle

\begin{abstract}
The \emph{wisdom of crowds} is an umbrella term for phenomena suggesting that the collective judgment or decision of a large group can be more accurate than the individual judgments or decisions of the group members. A well-known example illustrating this concept is the competition at a country fair described by Sir Francis Galton, where the median value of the individual guesses about the weight of an ox resulted in an astonishingly accurate estimate of the actual weight. In spirit, this phenomenon resembles classical results in probability theory, such as the law of large numbers, and it relies on independent decision-making. The accuracy of the group’s final decision can be significantly reduced if, under the dynamics of social influence, it appears to be driven by a few influential agents.

In this paper, we consider a group of agents who initially possess uncorrelated and unbiased noisy measurements of a common state of the world, modeled as a scalar parameter. Assume these agents iteratively update their estimates according to a simple non-Bayesian learning rule, commonly known in mathematical sociology as the French-DeGroot dynamics or iterative opinion pooling. As a result of this iterative distributed averaging process, each agent arrives at an asymptotic estimate of the state of the world, with the variance of this estimate determined by the matrix of weights the agents assign to each other. Every agent aims at minimizing the variance of her asymptotic estimate of the state of the world; however, such variance is also influenced by the weights allocated by other agents. To achieve the best possible estimate, the agents must  then solve a game-theoretic, multi-objective optimization problem defined by the available sets of influence weights.

We consider a class of influence mechanisms in which agents interact over a strongly connected aperiodic directed graph. Similar to dynamical models of reflected appraisals, we assume that a static matrix of relative influence weights is predefined, and the actual weights that an agent assigns to others are obtained by scaling these constants by her self-confidence value. An agent aims to select her self-confidence value in a way that minimizes the variance of her asymptotic estimate of the state of the world. We characterize both the Pareto frontier and the set of Nash equilibria in the resulting game. Additionally, we examine asynchronous best-response dynamics for the group of agents and prove their convergence to the set of strict Nash equilibria.
\end{abstract}



\section{Introduction}

The wisdom of crowds~\citep{Surowiecki:2004} has been actively studied in the social and behavioral sciences. In its broadest sense, this term refers to phenomena suggesting that cooperative decision-making generally outperforms individual decision-making.
A typical mathematical formalization involves a group of agents collecting noisy measurements of the same quantity of interest. Basic inferential arguments suggest that there exist linear aggregative functions of such data that yield a noise variance reduction with respect to the  individual measurements.
A paradigmatic example is the weight guessing game at a county fair, described by~\cite{Galton1907} who observed that the median of the villagers' individual guesses about the weight of an ox produced a surprisingly precise estimate of the actual weight. Similar phenomena of collective intelligence have recently been reported in animal groups~\citep{Berdahl2013}.

However, the wisdom-of-crowds phenomena can be undermined by mechanisms of opinion aggregation in influence systems. This is prominently illustrated by the influential work of~\cite{GolubJackson:2010}, which demonstrates that a group's collective wisdom may vanish as a result of ``naive'' non-Bayesian learning, known as the French-DeGroot model of iterative opinion pooling and often used to model rational consensus~\citep{French:1956,DeGroot,Lehrer:1976}.
This occurs, for instance, when all agents have equally accurate estimations but the influence system is biased toward certain agents, resulting in a non-uniform aggregation. More generally, this can happen when estimations vary in quality and the influence system fails to properly incorporate this information into the aggregation mechanism, as discussed in~\cite{GolubJackson:2010}, \cite{Bullo.Fagnani.Franci:2020},  and the references therein. It should be noted that the dynamics of iterative distributed averaging prove to be more than just a theoretical construct; their predictive power has been demonstrated in experiments with both small social groups~\citep{FriedkinJiaBullo:2016,FriedkinBullo:17,friedkin2019mathematical,FRIEDKINPROBULLO2021} and social media~\citep{Kozitsin:2022,Kozitsin:2023}.

In this paper, similarly to \cite{GolubJackson:2010}, we assume that a finite set of agents with different estimation capabilities make uncorrelated unbiased noisy measurements of a common state of the world, modeled as a scalar value. Afterwards, the agents engage on an iterative non-Bayesian learning process on a social network and aggregate their initial estimates through a distributed averaging mechanism {\it \`a la} French-DeGroot. Specifically, we assume that the agents  keep on updating their opinion --- representing their current estimate of the unobserved state of the world --- to a convex combination of their own current opinion and of a weighted average of those of her neighbors.

Similarly to models of reflected appraisals \citep{FriedkinBullo:2014}, we assume that the matrix of relative interpersonal weights modeling the social network is fixed, representing the proportions in which agents allocate influence to each other. Our focus, however is on the effect of the choice of the agents' self-confidence values, as measured by the weight $z_i$ that each agent $i$ assigns to its own opinion in a convex combination with the weighted average of her neighbors' opinions. Agents with $z_i = 1$ are called stubborn, as they ignore others entirely, while those with $z_i \in [0,1)$ are regular, assigning positive weight $1 - z_i$ to their neighbors' opinions~\citep{Acemoglu.ea:2013}.

Assuming the matrix of relative interpersonal weights is irreducible and aperiodic, standard distributed averaging techniques imply that, for any self-confidence profile $z$, agents' opinions converge to random variables whose expected value equals the true state of the world. The inverse variance $u_i(z) = 1/v_i(z)$ of agent $i$'s asymptotic opinion quantifies the precision or accuracy of her final estimate. Notably, for fixed interpersonal weights and observation variances, this precision depends not only on $z_i$ but also on the self-confidence of all other agents.

We model the agents' choice of their self-confidence value as a strategic decision making process. Specifically, we study the strategic game with player set coinciding with the set of agents,  whereby every player $i$
selects an action $z_i \in [0,1]$, representing her self-confidence, to maximize her utility $u_i(z)$, being
the precision of her final estimate of the state of the world. For this game, we provide a complete characterization of the set of pure strategy (strict) Nash equilibria. We then  study the asymptotic behavior of the asynchronous best response dynamics~\citep{Blume:1995}, a strategic dynamic learning process whereby at each time a randomly selected agent updates her self-confidence value to the one that maximizes her utility against the current self-confidence profile of the other agents. A peculiarity of the strategic game in question arises from discontinuities in the utility functions $u_i(z)$ as any $z_j$ approaches $1$, i.e., when a regular agent becomes stubborn, as this event disconnects the underlying graph and dramatically changes the nature of the final outcome of the distributed averaging process. This creates significant technical challenges in the analysis of the  Nash equilibria of the game and the asymptotic behavior of the asynchronous best dynamics.

It is worth noting that the game studied in this paper is not a network game in the sense of \cite{Galeotti.ea:2010} and \cite[Chapter 9]{Jackson:2008}. More precisely, while it can be trivially considered a graphical game on the complete graph, it is not a graphical game on the graph describing the interconnection pattern of the social network. This is because the utility $u_i(z)$ of player $i$ depends not just on her own action and on the actions of her neighbors in the social network, but also on the actions of every other player. In other words, the precision of one agent's asymptotic estimate of the state of the world is influenced by the self-confidence values of all the agents, not just her own or her neighbors'. In fact, the game considered in this paper is to be considered as a particular network formation game (see., e.g., \cite[Chapter 11]{Jackson:2008}, \cite{Catalano.ea:2024}, and references therein), since the agents choices of their self-confidence values determine the row-stochastic matrix $W(z)$ that drives the distributed averaging process. Within the class of network formation games,
the one considered in this paper stands out due to the specific choice
of the utility function $u_i(z)$ that is the precision of one agent's final estimate of the state of the world. This choice fundamentally affects the set of outcomes (e.g., the Nash equilibria) and contrasts with alternative utility choices considered in the literature.

\subsection*{Contributions of this paper}

The contribution of this paper is threefold. First, we introduce a novel network formation game whereby every player $i$, corresponding to a node, chooses a
self-confidence parameter $z_i\in[0,1]$ to maximize the precision $u_i(z)$ of her final opinion in a French-DeGroot distributed averaging process. The process starts from uncorrelated random initial opinions with a common expected value and evolves under the row-stochastic matrix $W(z)$, whose entries are defined as $W_{ij}(z) = (1 - z_i) P_{ij}$ for $i \ne j$, and $W_{ii}(z) = z_i + (1 - z_i) P_{ii}$. Here, $P$ is a fixed irreducible, aperiodic row-stochastic matrix.

Second, we provide a detailed analysis of pure strategy Nash equilibria for this game. Our main results in this part are presented in Theorems~\ref{theo:internalNash} and~\ref{theo:other-Nash}, as well as in Corollary~\ref{coro:no-nonstrct}. Essentially, our analysis shows that the considered game always admits pure strategy Nash equilibria corresponding to self-confidence choices that do not disconnect the graph (that is, no player chooses the maximal self-influence weight). All these equilibria are strict and equivalent in that, under each of them, every agent achieves the same optimal estimation -- the one they would obtain through optimal direct cooperation. For generic variance values, no other Nash equilibria exist. In contrast, when subsets of at least two agents share the same variance, additional non-strict Nash equilibria may emerge, which emerge from disconnecting a group of agents with identical variances from the rest.

Our third contribution is the analysis of a class of strategic learning dynamics for self-confidence adaptation. Unlike our conference paper~\cite{ComoFagnaniPro:2022}, we consider discrete-time  asynchronous best-response dynamics and prove their convergence to the set of  strict Nash equilibria for any initial condition (Theorem~\ref{theo:convergence}).

\subsection*{Structure of the paper}

We conclude with a brief outline of this work. In the final part of this introduction, we gather some notational conventions and graph-theoretic notions used throughout the paper. Section~\ref{sec:start} presents the model, in particular the estimation-driven network formation game.  In Section~\ref{sec:pareto}, we  provide preliminary results on the Pareto frontier of the game. Section \ref{sec:main} is the main technical section, in which we develop a fundamental analysis of the best response correspondences of this game and characterize the pure strategy Nash equilibria.
In Section~\ref{sec:dynam}, we study learning dynamics for self-confidence adaptation and prove its convergence to the set of strict Nash equilibria that determine the social optimum. We also include a brief numerical simulation to support the theoretical results. A conclusion section completes the paper.

\subsection*{Notational convections and graph-theoretic notions}\label{sec:notation}

Unless otherwise stated, all vectors are considered as columns.  We use $\mathbf{0}$ and $\one$ to denote the column vectors of all zeros and all ones, respectively. The symbol $\delta^i$ denotes the column of Kronecker symbols with $\delta^i_j=1$ when $j=i$ and $\delta^i_j=0$ otherwise; the dimensions of both vectors will be clear from the context.
For two vectors $x$ and $y$ in $\mathbb{R}^n$, the inequalities $\leq$,$<$,$\geq$, and $>$ are meant to hold entry-wise.
For a vector $z$ in $\mathbb{R}^n$, the symbol $[z]$ denotes the diagonal matrix with diagonal entries $z_1,\ldots,z_n$.
Throughout the paper, the transpose of a matrix $M$ is denoted by $M'$. As usual, a \emph{row-stochastic} matrix is a nonnegative square matrix $M$ whose rows sum up to $1$ i.e., such that $M\one=\one$.

A (finite, directed) \emph{graph} is the pair $\mc G=(\mc V,\mc E)$, of a finite set of nodes $\mc V$ and a set of (directed) links $\mc E\subseteq\mc V\times\mc V$.
Undirected graphs are regarded as special cases of graphs in which $(i,j)\in\mc E$ if and only if $(j,i)\in\mc E$.
The \emph{in-degree} and \emph{out-degree} of a node $i$ in a graph $\mc G$ are defined as
$d_i^-=|\{j\in\mc V:\,(j,i)\in\mc E\}|$  and  $d_i^+=|\{j\in\mc V:\,(j,i)\in\mc E\}|$, respectively. For $k\ge0$, graph $\mc G=(\mc V,\mc E)$ is $k$-regular if $d_i^-=d_i^+=k$,  for every $i$ in $\mc V$.

A length-$l$ \emph{walk} from a node $i$ to a node $j$ in a graph $\mc G=(\mc V,\mc E)$ is an $(l+1)$-tuple of nodes $(\gamma_0,\gamma_1,\ldots,\gamma_l)$ such that $\gamma_0=i$, $\gamma_l=j$, and $(\gamma_{k-1},\gamma_{k})\in\mc E$ for  $1\le k\le l$. A walk $(\gamma_0,\gamma_1,\ldots,\gamma_l)$ is referred to as a \emph{path} if $\gamma_h\ne\gamma_k$ for every $0\le h<k\le l$, except for possibly $\gamma_0=\gamma_l$, in which case the path is referred to as a \emph{cycle}.

For a graph $\mc G=(\mc V,\mc E)$ and a subset $\mc S\subseteq\mc V$, we let $\mc G[\mc S]=(\mc S, \mc E[\mc S])$ be the graph with node set $\mc S$ such that, for every $i\ne j$ in $\mc S$, $(i,j)\in\mc E[\mc S]$ if and only if there exists a path from $i$ to $j$ in $\mc G$ that does not pass through any intermediate node $k$ in $\mc S$ (in particular, if $(i,j)\in\mc E$ and $i,j\in\mc S$, then $(i,j)\in\mc E[\mc S]$).

A graph $\mc G=(\mc V,\mc E)$ is referred to as:  \emph{strongly connected} if for every two nodes $i\ne j$ in $\mc V$ there exists a path from $i$ to $j$; \emph{aperiodic} if the greatest common divisor of the lengths of all its cycles is equal to $1$.
A \emph{directed ring} is a $1$-regular, strongly connected graph $\mc G=(\mc V,\mc E)$.

Finally, to a nonnegative square $n\times n$ matrix $M$ we can associate the graph $\mc G_M=(\mc V,\mc E)$ with node set $\mc V=\{1,\ldots,n\}$ and the link set
$\mc E=\{(i,j):M_{ij}>0\}$. The matrix $M$ is \emph{irreducible} if its associated graph $\mc G_M$ is strongly connected, and \emph{aperiodic} if $\mc G_M$ is aperiodic.

\section{Problem setup and preliminary results}\label{sec:start}

Consider a finite set of agents $\calV=\{1,\ldots,n\}$, interacting over a social network. The latter is defined by a fixed row-stochastic matrix $P$ in $\mathbb{R}^{n\times n}$ of \emph{relative interpersonal weights}, representing the proportions in which agents allocate influence to each other. The social network's topology is then identified by the graph  $\mc G_P$.
Every agent is characterized by a real scalar \emph{opinion}  $x_i$ in $\R$ and a self-confidence value $z_i$ in $[0,1]$.
Upon stacking all the agents' opinions and self-confidence values respectively into the vectors $x$ (the opinion profile) in $\R^n$ and $z$ (the self-confidence profile) in the hypercube $$\mc Z\triangleq [0,1]^n.$$

\subsection{Wisdom of Agents and Distributed Estimation}

As in \cite{GolubJackson:2010}, we  assume that the initial opinions are the agents' guesses about some unknown common state of the world, modeled as a scalar parameter $\theta$ in $\R$. More precisely, let
\be\label{initial-opinions}x_i(0)=\theta+\xi_i,\qquad i\in\mc V\,,\ee
where $\xi_i$, for all $i$, are zero-mean uncorrelated random noise variables with positive finite variance $\sigma_i^2$, i.e.,
\be\label{noise-assumption}\E[\xi]=0\,,\qquad \E[\xi\xi']=[\sigma^2]\,,\ee
where $\sigma^2=(\sigma_1^2,\ldots,\sigma_n^2)$ denotes the column vector of the agents' variances. The inverse value $\sigma_i^{-2}$ of the variance or agent $i$, commonly referred to as her precision,  may be interpreted as a measure of her ``expertise'': the larger is the precision $\sigma_i^{-2}$ of agent $i$ (i.e., the smaller is her variance $\sigma_i^2$), the closer (in the mean square sense) is her initial opinion is to the actual state of the world $\theta$.


It is then assumed that the dynamics of opinions obey
the standard French-DeGroot model, that is, the following discrete-time system
\beq\label{eq.degroot}
x(t+1)=W(z)x(t),\;\;\qquad t=0,1,\ldots\,,\ee
where $W(z)$ is the stochastic matrix
\be\label{eq:W(z)}
W(z)\triangleq (I-[z])P+[z]\,.
\eeq
The above recursion requires every agent $i$ in $\mc V$ to update her opinion to a convex combination
$$x_i(t+1)=z_ix_i(t)+(1-z_i)\sum_{j\in\mc V}P_{ij}x_j(t)\,,$$
of her own current opinion and on a weighted average of those of her neighbors. Specifically, agent $i$'s self-confidence value $z_i$ represents the weight she places on her own current opinion $x_i(t)$, while the complementary influence weight $1 - z_i$ is distributed among her neighbors' current opinions $x_j(t)$ in the proportion determined by the relative interpersonal weights $P_{ij}$. We refer to an agent $i$ as \emph{stubborn} if $z_i=1$; in this case, the stubborn agent keeps her opinion constant in time, i.e., $x_i(t)=x_i(0)$ for every $t=0,1,\ldots$.
Let
\be
\mc S(z)\triangleq\{i=1,\ldots,n:\,z_i=1\}\,,\ee
stand for the set of all stubborn agents under the self-confidence profile $z$.

Throughout, we shall restrict our analysis to social networks whose influence matrix $P$ is irreducible, i.e., it is such that the associated graph $\mc G_P$ is strongly connected.
It follows from classical Perron-Frobenius theory~\cite[Ch.2, Th.~1.3]{BermanPlemmons_Book} that in this case there is a unique left eigenvector $\pi$ in $\R_+^{n}$ for $P$, which is a probability distribution, i.e.,
\be\label{pi}\pi=P'\pi\,,\qquad \one'\pi=1\,,\ee
and, furthermore, we have that\be\label{pi>0}\pi>0\,.\ee We refer to $\pi$ as the (weighted eigenvalue) \emph{centrality vector} of the social network. This notion of centrality naturally arises as the stationary distribution of a Markov chain, serves as a measure of social power~\citep{French:1956,ProTempo:2017-1}, and has numerous other applications~\citep{FrascaIshiiTempo:2015,Como.Fagnani:2015}.
The asymptotic behavior of the opinion dynamics \eqref{eq.degroot} can then be characterized as follows. 
\begin{lemma}\label{lemma:limit}
Consider a social network with irreducible aperiodic influence matrix $P$ and let $\pi$ be its  centrality vector.
Then, for every self-confidence vector $z$ in $\mc Z$,
\be\label{limWt}\lim_{t\to+\infty}W(z)^t=H(z)\,,\ee
where $H(z)$ in $\R^{n\times n}$ is a stochastic matrix such that:
\begin{enumerate}
\item[(i)] if $\mc S(z)=\emptyset$, then \be\label{Hz}H(z)=\one p(z)'\,,\ee
has all rows equal to $p(z)'$, where
\be\label{def:p(z)}
 p(z)\triangleq\frac1{\gamma(z)}(I-[z])^{-1}\pi\,,\ee
and $\gamma(z)$ is the positive scalar defined by
\be\label{def:gamma(z)}\gamma(z)\triangleq\sum_{i\in\mc V}\frac{\pi_i}{1-z_i}\,;\ee
\item[(ii)] if $\mc S(z)\ne\emptyset$, then $H(z)$ is the unique solution of the linear system
\be\label{lin-sys}
\ba{rclcl}
\ds H_{ij}(z)&=&\ds\sum_{k\in\mc V}P_{ik}H_{ik}(z)&\ &\forall i\!\not\in\!\mc S(z)\,,\ j\!\in\!\mc V\,,\\[15pt]
\ds H_{ij}(z)&=&\delta^i_j &\ &\forall i\!\in\! \mc S(z)\,,\  j\!\in\!\mc V\,.
\ea
\ee
\end{enumerate}
\end{lemma}
\begin{IEEEproof}  See Appendix \ref{sec:proof-lemma:limit}. \qed\end{IEEEproof}

\begin{remark}\label{rem:power}
When $\mc S(z)=\emptyset$, the vector $p(z)$ defined in~\eqref{def:p(z)} can be interpreted as the vector of agents’ \emph{social powers}~\citep{ProTempo:2017-1}. Each entry $p_i(z)$ measures how agent~$i$'s initial opinion contributes to the group's final opinion
\[
x_1(\infty)=\ldots=x_n(\infty)=\sum_{i\in\mc V} p_i(z)\,x_i(0)=p(z)'x(0).
\]
By changing the self-confidence profile $z$, the agents thus adjust their social powers in accordance with~\eqref{def:p(z)} and~\eqref{def:gamma(z)}. The eigenvalue centrality from~\eqref{pi},  $\pi = p(0)$, corresponds to the case in which agents have no self-confidence and rely solely on
their neighbors' opinions.
\end{remark}

\begin{remark}\label{rem:centrality} Notice that, if $\mc S(z)\ne\emptyset$, then Lemma \ref{lemma:limit}(ii) implies that $H(z)$ depends only on $\mc S(z)$, i.e., for every  two self-confidence profiles $z^{(1)}$ and $z^{(2)}$ in $\mc Z$ we have
\be\label{onlySz}\mc S(z^{(1)})=\mc S(z^{(2)})\ne\emptyset\quad \Rightarrow\quad  H(z^{(1)})=H(z^{(2)})\,.\ee
In particular, if $\mc S(z)=\{i\}$, then $H(z)=\one(\delta^i)'$.
\end{remark}

\begin{remark}\label{rem:hitting} The following probabilistic interpretation is insightful.
Consider a discrete-time Markov chain with finite state space $\mc V$ and transition probability matrix  $W(z)=(I-[z])P+[z]$. It is known that, if $\mc S(z)=\emptyset$, then such a Markov chain converges in probability to a random variable distributed in accordance with the vector $$\frac1{\gamma(z)}(I-[z])^{-1}\pi\,,$$
where $\gamma(z)$ is defined as in \eqref{def:gamma(z)}.
On the other hand, if $\mc S(z)\ne\emptyset$, then the Markov chain 
gets absorbed in finite time in one of the states in $\mc S(z)$: in this case, $H_{ij}(z)$ coincides with the probability of absorption of the Markov chain in state $j$ when starting in state $i$.
\end{remark}

Lemma \ref{lemma:limit} has the following consequence on the asymptotic behavior of the French-DeGroot opinion dynamics.

\begin{proposition}\label{prop:convergece}
Consider a social network with irreducible aperiodic influence matrix $P$. Let $\pi$ be its centrality vector and let $z$ in $\mc Z$ be a self-confidence vector. Assume that the initial opinion profile has entries as in \eqref{initial-opinions} for some $\theta$ in $\R$ and a random vector $\xi$ satisfying assumption \eqref{noise-assumption}.
Then, the French-DeGroot opinion dynamics \eqref{eq.degroot} is such that
\be\label{limxi}\lim_{t\to+\infty}x_i(t)=\theta+\zeta_i\,,\ee
for every agent $i$ in $\mc V$, where $\zeta_i$ are random variables with expected value
\be\label{Ezeta=0}\E[\zeta_i]=0\,,\ee
and variance
\be\label{estimation-variance}\ups_i(z)=\E[\zeta_i^2]=\sum_{j\in\mc V}H_{ij}^2(z)\sigma_j^2\,,\ee
where $H(z)$ in $\R^{n\times n}$ is the stochastic matrix from~\eqref{limWt}.
\end{proposition}
\begin{IEEEproof}
It follows from Lemma \ref{lemma:limit} that
$$x(t)=W(z)^tx(0)\xrightarrow[t\to+\infty]{}H(z)x(0)\,. $$
It follows from assumption \eqref{initial-opinions} and the fact that $H(z)$ is a row-stochastic matrix that
$$H(z)x(0)=\theta H(z)\one+H(z)\xi=\theta\one+\zeta\,,$$
where $\zeta=H(z)\xi$.
Thanks to assumption \eqref{noise-assumption}, we get that the entries of $\zeta_i$ are zero-mean, thus proving equation \eqref{Ezeta=0}, and have variance
$$\E[\zeta_i^2]=\sum_{j\in\mc V}\sum_{k\in\mc V}H_{ij}(z)\E[\xi_j\xi_k]H_{ik}(z)=\sum_{j\in\mc V}H_{ij}^2(z)\sigma_j^2\,,$$
thus proving equation \eqref{estimation-variance}.\qed
\end{IEEEproof}


%
%

The zero-mean random variable $\zeta_i$ in Proposition \ref{prop:convergece} represents the remaining noise in agent  $i$'s estimation of the state of the world $\theta$ at the end of the social interaction: its variance is given by \eqref{estimation-variance}.
We shall interpret \be\label{estimation-cost}\ups_i(z)=\sum_{j\in\mc V}H_{ij}^2(z)\sigma_j^2\,.\ee as the asymptotic (mean square) estimation cost of agent~$i$. Given the influence matrix $P$ and the agents' initial variance values $\sigma_i^2$, such estimation cost depends only on the self-confidence vector $z$. Henceforth, we adopt a common abuse of notation by identifying the vector $z$ in $\mc Z$ with the pair $(z_i, z_{-i})$ and writing
$$\ups_i(z)=\ups_i(z_i,z_{-i})\,,$$
to emphasize that the asymptotic estimation cost depends both on the agent's own self-confidence value $z_i$ and on the self-confidence profile of the other agents $$z_{-i}=(z_j)_{j\ne i}\in [0,1]^{n-1}.$$

\subsection{Problem Setup}

Assume that the agents are rational entities who select their self-confidence values $z_i$ to optimize the accuracy of their estimation of true state of the world $\theta$. Formally, we consider the strategic game with player set $\mc V$ and, for every player $i$ in $\mc V$, action set $[0,1]$ and utility function $$u_i(z)=1/v_i(z)\,,$$ coinciding with the inverse of the variance,  i.e., with the precision, of player $i$'s asymptotic estimate of the state of the world.

We shall use the following multi-objective optimization and  game-theoretic notions.
\begin{definition}
A self-confidence profile $z$ in $\mc Z$ is
\begin{enumerate}
\item[(i)] \emph{Pareto-optimal} if there exists no other self-confidence profile $\ov z$ in $\mc Z$ such that
$$\ups_i(\ov z)\le\ups_i(z)\,,\qquad\forall i\in\mc V\,,$$ and there exists an agent $j$ in $\mc V$ such that $$\ups_j(\ov z)<\ups_j(z)\,;$$
\item[(ii)] a \emph{Nash equilibrium} if for every agent $i$ in $\mc V$
$$\ups_i(z_i,z_{-i})\le\ups_i(\ov z_i,z_{-i})\,,\qquad \forall\ov z_i\in[0,1]\,;$$
\item[(iii)] a \emph{strict Nash equilibrium} if for every agent $i$ in $\mc V$
$$\ups_i(z_i,z_{-i})<\ups_i(\ov z_i,z_{-i})\,,\qquad \forall\ov z_i\in[0,1]\setminus\{z_i\}\,.$$
\end{enumerate}
We shall denote by $\mc P\subseteq\mc Z$ the Pareto frontier, i.e., the set of Pareto optimal self-confidence profiles and by $\mc N$ and $\mc N^*$, respectively, the set of Nash equilibria and the set of strict Nash equilibria.
\end{definition}

\textbf{Problem statement.} The goals of the following sections are twofold. We first aim to provide an exhaustive characterization of the Pareto frontier
(Section~\ref{sec:pareto}) and the structure of the Nash equilibria (Section~\ref{sec:main}). Then, we examine the convergence of asynchronous best-response dynamics in the game (Section~\ref{sec:dynam}).

\section{Pareto-optimal Self-confidence Profiles}\label{sec:pareto}
\color{black}
Proposition~\ref{prop:convergece} and Lemma~\ref{lemma:limit}(i) imply that, when all the agents $i$ choose self-confidence value $z_i<1$, they all share the same cost function
\be\label{unique-cost}\ups_i(z)=V(z)=\left(\sum\limits_{j\in\mc V}\frac{\pi_j}{1-z_j}\right)^{-2}\sum\limits_{j\in\mc V}\frac{\pi_j^2\sigma_j^2}{(1-z_j)^2}\,.\ee

Furthermore, one can find the self-confidence profiles $z$ minimizing the common cost $V(z)$ as follows.
\begin{lemma}\label{lemma:min}
If  $\sigma_i^2>0\,$ for every $i=1,\ldots,n$, then
\be\label{eq:Vmin}
\min_{\substack{\mu\in\R_+^n:\\\one'\mu=1}}\sum_{i\in\mc V}\mu_i^2\sigma_i^2=V_{min}:=\left(\sum_{i\in\mc V}\sigma_i^{-2}\right)^{-1} \,,
\ee
and the minimum above is achieved in the sole probability distribution $\mu^*$ with entries
\be\label{mu}\mu^*_i=\frac{\sigma_i^{-2}}{\ds\sum_{j\in\mc V}\sigma_j^{-2}}\,,\ee
for every $i=1,\ldots,n$.
\end{lemma}
\begin{IEEEproof} See Appendix~\ref{sec:proof-lemma:min}.\footnote{An alternative proof earlier appeared in~\cite{Peluffo:2019}.} 
\end{IEEEproof}

A hypothetical centralized social planner who has full access to all agents' measurements \eqref{initial-opinions} would choose the unbiased linear estimator $\sum_{i\in\mc V}\mu_i^*x_i(0)$
in order to minimize error variance. As a consequence, every vector $z$ in $\mc Z$ such that $H(z)=\one(\mu^*)'$ is Pareto optimal, and
if such vectors exist, they exhaust the Pareto frontier $\mc P$. This is formalized in the following result.
\begin{proposition}\label{prop:Pareto}
Consider a social network with irreducible aperiodic influence matrix $P$ and centrality vector $\pi$.
Assume that the agents make uncorrelated initial estimations with variance $\sigma_i^2>0$ for every $i$ in $\mc V$.
Let
\be
\label{Z*}\mc Z^*\triangleq\{\one-\alpha[\pi]\sigma^2:\,0<\alpha\leq\alpha_*\}\,,\ee
where
\be\label{alpha*}\alpha_*\triangleq 1/{\max_i\{\pi_i\sigma_i^2\}}\,.\ee
Then, $$\mc P=\mc Z^*\subseteq\mc N^*\,.$$
Moreover, for every Pareto-optimal vector $z$ in $\mc Z^*$, the social powers associated to the French-DeGroot model~\eqref{eq.degroot} (see Remark~\ref{rem:power}) are $p_i(z)=\mu_i^*$, and the variance of the final consensus opinion is found as $V_{min}$, defined in~\eqref{eq:Vmin}.
\end{proposition}
\begin{IEEEproof}
We first notice that when stubborn agent are present, i.e., when $\mc S(z)\ne\emptyset$, the matrix $H(z)$ cannot be of the type $\one(\mu^*)'$,
as in this case $H(z)$ has zero entries by equation \eqref{lin-sys} of Lemma \ref{lemma:limit}.
Instead, if we consider the expression of $H(z)$ in \eqref{Hz} for the case when there are no stubborn agents (i.e., when $z<\one$), we see that
$$H(z)=\one(\mu^*)'\quad\Leftrightarrow\quad \gamma(z)^{-1}(I-[z])^{-1}\pi=\mu^*$$
Because of the way $\mu^*$ is defined (see \eqref{mu}), the right hand side equality above holds true if and only if the vector $z$ satisfies the relations
\be\label{Pareto}1-z_i=\alpha\pi_i\sigma_i^2\,,\ee
for every $i$ in $\mc V$, for some constant $\alpha$. Notice that \eqref{Pareto} has a solution in $[0,1)^n$ if and only if $0<\alpha\leq\alpha_*$.
Hence, $H(z)=\one(\mu^*)'$ if and only if $z\in\mc Z^*$ (in particular, the set of such self-confidence profiles $z$ is nonempty). As has been already noted, this entails that $\mc P=\mc Z^*$. On the other hand, it directly follows from the definition of $\mc Z^*$ in \eqref{Z*} that for every $z$ in $\mc Z^*$, any unilateral modification of the self-confidence by one agent would lead to a vector
$\bar z\not\in \mc Z^*$. Hence $z$ is a strict Nash equilibrium.\qed
\end{IEEEproof}

\begin{remark} Proposition~\ref{prop:Pareto} shows that the Pareto frontier $\mc P$ is a segment within $[0,1)^n$. All vectors in $\mc P$ are equivalent from a cost perspective: in such configurations, every player achieves the same cost, corresponding to the minimum of the convex program presented in Lemma~\ref{lemma:min}.
\end{remark}

\section{Best-Response Maps and Nash Equilibria}\label{sec:main}

To refine our analysis of the game, we first undertake a fundamental analysis of best response correspondences.

\subsection{The best response correspondence and its properties}
We denote by $\mathcal B_i(z_{-i})\subseteq [0,1]$ the best response (BR) set of an agent $i$ in $\mc V$ when the rest of the agents have chosen the self-confidence profile $z_{-i}$, i.e.,
\[
\mathcal B_i(z_{-i})\triangleq\argmin_{z_i\in[0,1]}\ups_i(z_i,z_{-i}).
\]

It is convenient to set some notation. We define the following aggregated expressions:
\beq\label{eq.aux0}
A(z)
=\sum_{j\in\mc V}\frac{\pi_j}{1-z_j},\qquad
B(z)
=\sum_{j\in\mc V}\frac{\pi_j^2\sigma_j^2}{(1-z_j)^2}.
\eeq
\beq\label{eq.aux1}
A_{{k}}(z_{-k})
=\sum_{j\ne k}\frac{\pi_j}{1-z_j},\qquad
B_{{k}}(z_{-k})
=\sum_{j\ne k}\frac{\pi_j^2\sigma_j^2}{(1-z_j)^2}.
\eeq
It is also convenient to introduce new variables
\beq\label{change}
y_i\dfb\frac{1}{1-z_i}\in [1,+\infty),\qquad i\in\mc V\,.
\eeq
For every function of the vector variable $z$, say $f(z)$, we use the convention of indicating with a bar on top the result   $\ov f(y)$ of the composition of the function  $f$ with the variable change {determined by the inverse of transformation} \eqref{change}.
So, in particular, we indicate with $\ov V(y)$ the common cost as defined in \eqref{unique-cost} and with
$\ov A (y), \ov B(y), \ov A_{{k}}(y_{-k}), \ov B_{{k}}(y_{-k})$ the expressions corresponding to \eqref{eq.aux0} and \eqref{eq.aux1} with this substitution.

For every agent $i$ in $\mc V$, the common cost can be written as
\[
\ov V(y)=\frac{\ov B(y)}{\ov A(y)^2}=
\frac{\pi_i^2\sigma_i^2y_i^2+\ov B_{{i}}(y_{-i})}{\left(\pi_iy_i+\ov A_{{i}}(y_{-i})\right)^2}.
\]
Differentiating $\ov V$ with respect to $y_i$, we obtain
\be\label{eq.partial}
\ba{rcl}
\ds\frac{\partial\ov V}{\partial y_i}(y)
&=&\ds\frac{2\pi_i^2\sigma_i^2y_i}{(\pi_iy_i+\ov A_{i}(y_{-i}))^2}-
\frac{2(\pi_i^2\sigma_i^2y_i^2+\ov B_{i}(y_{-i}))\pi_i}{(\pi_iy_i+\ov A_{i}(y_{-i}))^3}
\\[10pt]
&=&\ds\frac{2\pi_i\left(\ov A_{i}(y_{-i})\pi_i\sigma_i^2y_i-\ov B_{i}(y_{-i})\right)}{\left(\pi_iy_i+\ov A_{i}(y_{-i})\right)^3}.
\ea
\ee

Our study begins with the case in which all agents $j$ in $\mc V\setminus\{i\}$ choose a self-confidence value $z_j < 1$.
The following result determines the form of agent~$i$’s best response in this case and shows that if no other agent $j$ in $\mc V \setminus \{i\}$ is stubborn, then it is never a best response for agent $i$ to be stubborn.
\begin{proposition}\label{prop:brm}
The best response of every agent $i$ in $\mc V$ to a self-confidence profile $z_{-i}$ such that $\mc S(z_{-i})=\emptyset$ is given by
\beq\label{eq.brm}
\mc B_i(z_{-i})=\left\{\left[1-\frac{ A_{i}(z_{-i})\pi_i\sigma_i^2}{B_{i}(z_{-i})}\right]_+\right\}\,.
\eeq
{In particular, $\mc B_i(z_{-i})\subseteq [0, 1)$.}
\end{proposition}
\begin{IEEEproof}
Observe that the assumption $\mc S(z_{-i})=\emptyset$ is equivalent to that $0\le z_{-i}<\one$.
We first analyze $\ups_i(z_i,z_{-i})$ when $z_i<1$. In this case, $\ups_i(z_i,z_{-i})=V(z)$, where  $V(z)$ is the common cost function defined in
\eqref{unique-cost} that we analyze through the change of variables \eqref{change}.
Consider the expression for the derivative of $\ov V$ in (\ref{eq.partial}).
For a given $y_{-i}$, the function $y_i\mapsto\ov V(y_i,y_{-i})$ is strictly decreasing when $1\leq y\le y_i^*$ and strictly increasing when $y\ge y_i^*$, where
$$y_i^*=\min\left\{1,\frac {\bar B_i(y_{-i})}{\ov A_{{i}}(y_{-i})\pi_i\sigma_i^2}\right\}\,.$$
Hence, $y_i^*$ is the strict minimizer of $\bar\ups_i(y_i,y_{-i})$ on the interval $[1,+\infty)$.
This minimizer corresponds to
\[
z_i^*= 1-\frac{1}{y^*_i}=\left[1-\frac{ A_{{i}}(z_{-i})\pi_i\sigma_i^2}{B_i(z_{-i})}\right]_+\in{[0,1)}\,.
\]

Now, observe that, for  every fixed $0\le z_{-i}<\one$, Lemma \ref{lemma:limit}(i) implies that
$$\lim_{z_i\uparrow 1}H(z_i,z_{-i})=\one(\delta^i)'=H(1,z_{-i})\,,$$
so that,  by \eqref{estimation-cost}, the function  $z_i\mapsto\ups_i(z_i,z_{-i})$ is continuous on $[0,1]$. This implies that $z_i^*$ is in fact the unique minimizer of $z_i\mapsto\ups(z_i,z_{-i})$ over closed interval $[0,1]$.\hfill\qed
\end{IEEEproof}

In the case where the self-confidence profile $z_{-i}$ of the other agents contains some entries equal to $1$ (whose corresponding agents are stubborn), the analysis of the best response $\mc B_i(z_{-i})$ of agent $i$ is more involved. We start with the following result.

\begin{proposition}\label{prop.brm1}
For every agent $i$ in $\mc V$ and self-confidence profile $z_{-i}$ such that $\mc S(z_{-i})\neq\emptyset$, we have that:
\begin{enumerate}
\item[(i)] $z_i\mapsto\ups_i(z_i,z_{-i})$ is constant on $[0,1)$; \\
\item[(ii)] $\mc B_i(z_{-i})\in\{[0,1), [0,1], \{1\}\}$;\\
\item[(iii)] if $\sigma_i^2\ge\sigma_j^2$ for every $j$ in $\mc S(z_{-i})$, then $$\mc B_i(z_{-i}))\in\{[0,1), [0,1]\}\,;$$
\item[(iv)] if $\sigma_i^2>\sigma_j^2$  for every $j$ in $\mc S(z_{-i})$, then $$\mc B_i(z_{-i})=[0,1)\,.$$
\end{enumerate}
\end{proposition}
\begin{IEEEproof}
We consider the general expression \eqref{estimation-cost} for the utility of agent $i$, $v_i(z_i, z_{-i})$. As the entries of the matrix $H$ depend on $z$ only through the set $\mc S(z)$ (see Remark \ref{rem:centrality}), and since $\mc S(z)$ is constant for $z_i$ varying in $[0,1)$, item (i) is proven. Item (ii) is immediate from (i).

 In order to prove the remaining points, we consider the following estimation that holds true for every $z_i$ in $[0,1)$ and if $\sigma_i^2\ge\sigma_j^2$ for every $j$ in $\mc S(z_{-i})$.
\begin{equation}\label{eq:1-optimal}
\ba{rcl}
\ds\ups_i(z_i,z_{-i})&=&\ds\sum_{j\in \mc S(z_{-i})}H_{ij}^2(z)\sigma_j^2\\
&\leq&\ds\sum_{j\in \mc S(z_{-i})}H_{ij}(z)\sigma_j^2\\
&\leq&\ds\sigma_i^2\\
&=&\ds\ups_i(1,z_{-i})
\ea
\end{equation}
We notice that the first equality follows from equation~\eqref{estimation-cost} and the fact that $H_{ij}(z)=0$ for every agent $j\not\in\mc S(z)$, the first inequality follows from the fact that $H_{ij}(z)\in [0,1]$ for every $i,j$, and the second inequality from the assumption made on variances and the fact that $H(z)$ is a stochastic matrix.

Inequality~\eqref{eq:1-optimal} implies that any $z_i$ in $[0, 1)$ is necessarily a best response and thus yields item (iii). We further notice that for the second inequality in expression \eqref{eq:1-optimal} to be an equality, it must be that $\sigma_j^2=\sigma_i^2$ for every $j\in \mc S(z_{-i})$ such that $H_{ij}(z)>0$. Recalling that $H(z)$ is stochastic, such an index $j$ surely exists and this implies that under the assumption in item (iv), the second inequality in expression~\eqref{eq:1-optimal} is strict. This proves (iv).
\qed
\end{IEEEproof}

Proving statements (iii) and (iv) in Proposition~\ref{prop.brm1}, we have in fact established a stronger result.
\begin{proposition}\label{prop:reachable}
Let $i$ in $\mc V$ be an agent with minimal wisdom, i.e.,
$\sigma_i^2\geq\sigma_j^2$ for every $j$ in $\mc V$. Consider a self-confidence profile $z_{-i}$ such that
$\mc S(z_{-i})\neq\emptyset$.
Then, the following conditions are equivalent:
\begin{enumerate}
\item[(i)] $1\in \mc B_i(z_{-i})$
\item[(ii)] $[0,1]=\mc B_i(z_{-i})$
\item[(iii)] for every $z_i$ in $[0,1)$,
there exists only one node $j$ in $\mc S(z_{-i})$ reachable from node $i$ in the graph $\mc G_{W(z)}$ and, for this node $j$, one has $\sigma_i^2=\sigma_j^2$.
\end{enumerate}
\end{proposition}
\begin{IEEEproof} We prove (i)\,$\Rightarrow$\,(ii)\,$\Rightarrow$\,(iii)\,$\Rightarrow$\,(i)

(i)\,$\Rightarrow$\,(ii) follows from (iii) of Proposition \ref{prop.brm1}.

(ii)\,$\Rightarrow$\,(iii): under assumption (ii) we must have all equalities in expression \eqref{eq:1-optimal}. This implies that $H_{ij}(z)\in\{0,1\}$ for every $j\in \mc S(z_{-i})$ and that $\sigma_i^2=\sigma_j^2$ for every $j$ for which $H_{ij}(z)=1$. Since $H(z)$ is a stochastic matrix, this implies that there exists just one $j$ such that $H_{ij}(z)=1$ and for such $j$ it holds $\sigma_i^2=\sigma_j^2$. Considering the probabilistic interpretation of $H_{ij}(z)$ as the absorbing probability in the state $j$ of a Markov chain with matrix $W(z)$ initialized in $i$, it follows that
$H_{ij}(z)>0$ yields the existence of a path in $\mc G_{W(z)}$ from $i$ to $j$.

(iii)\,$\Rightarrow$\,(i): using the same probabilistic interpretation, we deduce from (iii) that $H_{ij}(z)=1$. This fact and the assumptions on the variances imply that all inequalities in expression \eqref{eq:1-optimal} turn into equalities. This yields (i).
\qed
\end{IEEEproof}

\subsection{The structure of Nash equilibria}

We are now ready to study the Nash equilibria of the game. First, we show that the only Nash equilibria in $[0,1)^n$ are strict and coincide with the Pareto frontier we have already examined.
 \begin{theorem}\label{theo:internalNash} The set $\mc P=\mc Z^*$ comprises all Nash equilibria in $[0,1)^n$, that is,
 $\mc N\cap [0,1)^{n}=\mc N^*\cap [0,1)^{n}= \mc Z^*$.
 \end{theorem}
 \begin{IEEEproof}
Notice first that the inclusion $\mc N^*\cap [0,1)^{n}\supseteq \mc Z^*$ follows from Proposition \ref{prop:Pareto}.
To complete the proof, it remains to show that there are no Nash equilibria in the set $[0,1)^{n}\setminus\mc Z^*$.

We first show that if $z=\mathbf{0}$ is a Nash equilibrium, then $\mathbf{0}\in\mc Z^*$. By the definition of a Nash equilibrium, $0\in\mc B_i(\mathbf{0})$ for all $i$, meaning that $z_i=0$ minimizes $V(z_i,\mathbf{0})$ and $y_i=1$ minimizes $\bar V(y_i,\one)$. In view of~\eqref{eq.partial}, one has
\[
\frac{\partial\bar V}{\partial y_i}(\one)=\frac{2\pi_i\left(\ov A(y)\pi_i\sigma_i^2y_i-\ov B(y)\right)}{\left(\pi_iy_i+\ov A_{i}({y_{-i}})\right)^3}\Big\vert_{y=\one}\geq 0.
\]
One notices, however, that the numerators of the latter fractions,
for $y=\one$, satisfy the following property:
$$\sum\limits_{i\in\mc V}2\pi_i(\ov A(\one)\pi_i\sigma_i^2-\ov B(\one))=\sum\limits_{i\in\mc V}2\pi_i^2\sigma_i^2-\sum\limits_{i\in\mc V}2\pi_i^2\sigma_i^2=0.
$$
Hence all fractions must vanish at $y=\one$, entailing that
$$\ov B(\one)-\ov A(\one)\pi_i\sigma_i^2=0,\quad\forall i\in\mc V,$$
which implies that $\pi_i=\sigma_i^{-2}/\sum_{j\in\mc V}\sigma_j^{-2}$, that is, $\pi=\mu^*$. By \eqref{Z*}, we have $\mathbf{0}\in\mc Z^*$ (this vector corresponds to $\alpha=\alpha^*$ in \eqref{Z*}).

Consider now a general Nash equilibrium $\bar z$ in $[0,1)^n$. We observe that $\bar z$ also constitutes a Nash equilibrium in the game with the more restricted strategy set $\mc Z_{\bar z}=\prod_{i=1}^n[\bar z_i,1)$. For every vector $z$ in $\mc Z_{\bar z}$, the matrix of the French-DeGroot model~\eqref{eq.degroot} can be rewritten as follows
\[
\begin{gathered}
W(z)=(I-[\hat z])W(\bar z)+[\hat z],\;\;\text{where}\\
\hat z\triangleq(I-[\bar z])^{-1}(z-\bar z)\in[0,1)^n.
\end{gathered}
\]
Hence, $\mathbf{0}$ is a Nash equilibrium in the game, defined by the matrix $\bar P=W(\bar z)$ with the corresponding left eigenvector
$\ov\pi=\gamma(\ov z)^{-1}(I-[\ov z])^{-1}\pi$, which, as has been shown, means that $\ov\pi=\mu^*$. In the proof of Proposition~\ref{prop:Pareto}
we have seen that this holds true if and only $\ov z\in\mc Z^*$.
\hfill\qed
\end{IEEEproof}

Theorem \ref{theo:internalNash} does not rule out the presence of Nash equilibria on the boundary of $[0,1)^n$. The following result identifies certain necessary conditions for such boundary equilibria to exist.
\begin{theorem}\label{theo:other-Nash}
Consider a social network with irreducible aperiodic influence matrix $P$ and centrality vector $\pi$.
Assume that the initial measurements of the agents are uncorrelated with variance $\sigma_i^2>0$ for $i$ in $\mc V$.
Then, every possible Nash equilibrium $z$ such that $\mc S(z)\neq\emptyset$
is non strict and is such that:
\begin{enumerate}
\item[(i)] $|\mc S(z)|\geq 2$;
\item[(ii)] $\sigma_i^2=\sigma_j^2$ for every $i$ and $j$ in $\mc S(z)$;
\item[(iii)] the graph $\mc G_P[\mc S(z)]$ is a directed ring\footnote{Notice that the undirected connected graph with two nodes is considered as a special case of a directed ring.}.
\end{enumerate}
\end{theorem}

\begin{IEEEproof}

{If $\mc S(z)=\{i\}$, Proposition~\ref{prop:brm} implies that $1\notin\mc B_i(z_{-i})$ and therefore $z$ is not a Nash. This shows that $|\mc S(z)|\ne 1$ and proves the necessity of condition (i).

Denote $$\mc S^*(z)\triangleq\argmax_{\ell\in\mc S(z)} \sigma_{\ell}^2\,.$$
If $\mc S^*(z)\neq \mc S(z)$, consider the shortest path $\gamma$ in $\mc G_P$, starting at some vertex from $\mc S^*(z)$ and leading to another vertex from $\mc S(z)\setminus\mc S^*(z)$ (such a path exists since the graph is strongly connected).
Denoting the starting and final vertices $i$ and $j$ respectively, the other vertices visited by the path $\gamma$ are, by construction, out of $\mathcal{S}(z)$.
Consequently, for every configuration $z'$ such that $z'_{-i}=z_{-i}$ and $z'_i<1$, $\gamma$ remains a path in the graph $\mc G_{W(z')}$. Since $\sigma_i^2>\sigma_j^2$, this implies that condition (iii) in Proposition \ref{prop:reachable} does not hold. This yields $1\not\in\mc B_i(z_{-i})$ contradicting the fact that $z$ is a Nash equilibrium.
Hence, condition (ii) is necessary: $\sigma_j^2=\sigma_*^2$ for every $j$ in $\mc S(z)$.

The proof that (iii) is necessary relies on similar arguments.
Recall first the definition of the restricted graph $\mc G_P[\mc S(z)]=(\mc S(z),\mc E[\mc S(z)])$ and notice that, since $\mc G_P$ is strongly connected, so is $\mc G_P[\mc S(z)]$.
%
%
%
%
In particular, each node has at least one out-neighbor. On the other hand, the existence of an edge $(i,j)$ in $\mc G_P[\mc S(z)]$ implies the existence of a path from $i$ to $j$ in $\mc G_P$ that does not visit other nodes from $\mc S(z)$. As this remains a path in $\mc G_{W(z')}$  for every configuration $z'$ such that $z'_{-i}=z_{-i}$ and $z'_i<1$, by virtue of (iii) in Proposition~\ref{prop:reachable}, we conclude that such node $j$ must be unique. In other terms, $i$ has just one out-neighbor. This implies condition (iii).

Finally, the fact that $z$ is not strict (as a Nash equilibrium) follows again from condition (ii) of Proposition~\ref{prop:reachable} applied to a stubborn agent $i$.
\qed

}
\end{IEEEproof}

Theorem \ref{theo:other-Nash} has the following direct consequences.
{\color{black}\begin{corollary}\label{coro:no-nonstrct}
The Pareto-optimal self-confidence profiles are the only strict Nash equilibria: $\mc N^*=\mc Z^*$. Additionally,
		if $\sigma_i^2\neq\sigma_j^2$ for $i\neq j$ in $\mc V$, then $$\mc N=\mc N^*=\mc Z^*$$
\end{corollary}
\begin{IEEEproof} By Theorem \ref{theo:other-Nash}, all strict Nash equilibria are inside $[0,1)^n$, in view of Theorem~\ref{theo:internalNash}, which proves the first statement. For mutually distinct $\sigma_i^2$, there are no other Nash equilibria due to statement (ii) in Theorem~\ref{theo:other-Nash}.
	\qed\end{IEEEproof} \medskip
}

 \begin{corollary}\label{coro:undirected} If $\mc G_P$ is undirected, then $|\mc S(z)|=2$ for every possible non-strict Nash equilibrium $z$.
\end{corollary}
\begin{IEEEproof}  Condition (iii) in Theorem \ref{theo:other-Nash} is never satisfied if $\mc G_P$ is undirected and $|\mc S(z)|\ge3$.\qed\end{IEEEproof}\medskip

While Theorem \ref{theo:other-Nash} and its corollaries provide necessary conditions for the existence of non-strict Nash equilibria {\color{black}that lay on the boundary of $[0,1)^n$}, the following result provides a sufficient condition.
\begin{proposition}\label{prop:main} Consider a social network with irreducible aperiodic influence matrix $P$ and assume that the initial measurements of the agents are uncorrelated with variance $\sigma_i^2>0$ for every $i$ in $\mc V$. Then, every self-confidence vector $z$ in $\mc Z$ such that
\begin{enumerate}
\item[(i)]  $|\mc S(z)|=2$,
\item[(ii)] 
$\sigma^2_j\le\sigma_i^2$  for every\footnote{This condition implies, in particular, that the agents from $\mc S(z)\triangleq\{k,\ell\}$ are equally wise:  $\sigma^2_k=\sigma^2_{\ell}$.} $j$ in $\mc S(z)$ and $i$ in $\mc V$,
\end{enumerate}
is a non-strict Nash equilibrium.
\end{proposition}
\begin{IEEEproof}
It follows from (iii) of Proposition \ref{prop.brm1} that
$\mc B_i(z_{-i})\supseteq[0,1)\ni z_i$ for every $i$ in $\mc V\setminus\mc S(z)$.
On the other hand, for every $i$ in $\mc S(z)$, it follows that condition (iii) in Proposition~\ref{prop:reachable} holds true. This entails  that
$\mc B_i(z_{-i})=[0,1]$ and thus surely $z_i\in \mc B_i(z_{-i})$. Thus, $z$ is a (non-strict) Nash equilibrium.
\qed
\end{IEEEproof}

\section{Asynchronous best response dynamics}\label{sec:dynam}
 In this section, we analyze the dynamic behavior of the asynchronous best response dynamics \cite{Blume:1995}.	
The asynchronous best response dynamics is a discrete-time Markov chain $Z(t)$,  with state space $[0,1]^{\mc V}$ and transition kernel described as follows.
At every discrete time step $t=0,1,2,\dots$,
an agent $k$ is chosen uniformly at random from the agent set $\mc V$ and she updates her action $Z_k(t)$ to a new action $Z_k(t+1)$ chosen uniformly at random from her best response set $\mc B_k(Z_{-k}(t))$, while the rest of the agents keep their actions unaltered, i.e., $Z_{-k}(t+1)=Z_k(t)$.

{Next lemma shows that, from any initial condition, the asynchronous best response dynamics $Z(t)$ gets inside $[0,1)^{\mc V}$ (no agent is stubborn) in finite time with probability one.}

\begin{lemma}\label{lemma:BR} For every distribution of the initial configuration $Z(0)$ on $[0,1]^{\mc V}$, there exists a nonnegative random time $T$ such that $\P(T<+\infty)=1$, and the asynchronous best response dynamics satisfies \be\label{Z<1}Z(t)<\one\,,\qquad \forall t\ge T\,.\ee
\end{lemma}
\begin{IEEEproof}
Consider the stochastic process $Q(t)$ on the finite state space $\mc Q=\{0,1\}^{\mc V}$, defined by
$$Q_i(t)=\lfloor Z_i(t)\rfloor\,,\qquad i\in\mc V\,,\ t=0,1,2,\ldots\,,$$
i.e., $Q_i(t)=1$ if agent $i$ is stubborn at time $t$ and $Q_i(t)=0$ otherwise. Hence, $$\mc S(Z(t))=\{i\in\mc V:\,Q_i(t)=1\}\,.$$

Recall from Proposition~\ref{prop:brm}
that the best response of every agent $i$ in $\mc V$
such that $Q_{-i}(t)=\0$, i.e.,  such that $\mc S(Z_{-i}(t))=\emptyset$, satisfies $\mc B_i(Z_{-i}(t))\subseteq [0,1)$.
This implies that, if the agent $i$ that gets randomly activated at time $t$ is such that $Q_{-i}(t)=\0$, then necessarily $Q_i(t+1)=0$.

On the other hand, for every every agent $i$ in $\mc V$ such that $Q_{-i}(t){\neq \0}$, i.e., such that $\mc S(Z_{-i}(t)){\neq}\emptyset$,
Proposition \ref{prop.brm1}(ii) ensures that $$\mc B_i(Z_{-i}(t))\in\{[0,1),[0,1],\{1\}\}\,.$$ Observe, that in this case, Equation \eqref{estimation-cost}, Lemma \ref{lemma:limit}(ii), and Remark \ref{rem:centrality} imply that
$$\ups_i(z_i,Z_{-i}(t))=\sum_{j\in\mc V}H_{ij}^2(z_i,Z_{-i}(t))\sigma_j^2\,,$$
depends on $Z_{-i}(t)$ only through $\mc S(Z_{-i}(t))$, i.e., through $Q_{-i}(t)$, hence so does $\mc B_i(Z_{-i}(t))$.
This implies that there exists a function $p_i:\{0,1\}^{\mc V\setminus\{i\}}\to[0,1]$, such that, if the agent $i$ that gets randomly activated at time $t$ is such that $Q_{-i}(t)\ne\0$, then  with conditional probability $p(Q_{-i}(t))$ we get $$\mc B_i(Z_{-i}(t))\in\{[0,1),[0,1]\}\,,$$  and hence $Z_i(t+1)=0$, while with conditional probability $1-p(Q_{-i}(t))$ we get $$\mc B_i(Z_{-i}(t))=\{1\}\,,$$ hence $Z_i(t+1)=1$.
In fact, Proposition \ref{prop.brm1}(iii) implies that, if
$$i\in\argmax\limits_{j: Q_j(t)=1}\sigma_j^2\,,$$
then $p(Q_{-i}(t))=1$.

The arguments above imply that $Q(t)$ is itself a Markov chain on $\mc Q$ such that
\be\label{absorbingQ}\P(Q(t+1)=\0|Q(t)=\0)=1\,,\ee
and that for every $q$ in $\mc Q\setminus\{\0\}$,
 there exists $q^-$ such that
 \be\label{q-1}\sum_{i\in\mc V}q^-_i=\sum_{i\in\mc V}q_i-1\,,\ee and
\be\label{Pq-1}\P(Q(t+1)=q^-|Q(t)=q)\ge\frac1n>0\,.\ee
It then follows from \eqref{absorbingQ}, \eqref{q-1}, and \eqref{Pq-1} that $\0$ is absorbing state for the Markov chain $Q(t)$ on the finite state space $\mc Q=\{0,1\}^{\mc V}$. Standard results on finite Markov chains~\citep[Chapter 1]{norris1998markov} imply that there exists a nonnegative random absorption time $T$ such that $\P(T<+\infty)=1$, and $Z(t)=\0$ for every $t\ge T$, so that
\eqref{Z<1} holds true.\qed \end{IEEEproof}

Finally, we prove the convergence of the asynchronous best response dynamics $Z(t)$ to the set $\mc Z^*$ of internal Nash equilibria.

\begin{theorem} \label{theo:convergence}
For every distribution of the initial configuration $Z(0)$ on $[0,1]^{\mc V}$, the sequence $Z(t)$ generated by the asynchronous best response dynamics converges to the set $\mc Z^*$, with probability one.
\end{theorem}
\begin{IEEEproof} We fix an initial condition $z(0)$  and we take $T$ as in Lemma \ref{lemma:BR}, so that $Z(t)\in [0,1)^n$ for every $t\geq T$.

We first show that every limit point of the sequence $Z(t)$ belongs to $[0,1)^n$.
To this aim, we define
$$M(t)=\max\limits_{j=1,\dots , n}\frac{\pi_j\sigma_j^2}{1-Z_j(t)}$$
and then we show that $M(t)$ is monotonically non-increasing for $t\geq T$. For a given time $t$, we pick
$$m\in \argmax\limits_{j=1,\dots , n}\frac{\pi_j\sigma_j^2}{1-Z_j(t)}\,,$$
and we notice that for any other player $k$ in $\mc V\setminus\{m\}$, the form of the best response in \eqref{eq.brm} implies that
$$\frac{\pi_k\sigma_k^2}{1-Z_k(t+1)}=\max\left\{ \pi_k\sigma_k^2, \frac{B_k(Z_{-k}(t))}{A_k(Z_{-k}(t))}\right\}\,.$$
Notice now that
$$\begin{array}{rcl}B_k(Z_{-k}(t))&=&\ds\sum\limits_{i\neq k}\frac{\pi_j^2\sigma_j^2}{(1-Z_j(t))^2}\\ &=&\ds\sum\limits_{i\neq k}\frac{\pi_j}{1-Z_j(t)}\frac{\pi_j\sigma_j^2}{1-Z_j(t)}\\ &\leq&\ds A_k(Z_{-k}(t))\frac{\pi_m\sigma_m^2}{1-Z_m(t)}\,.\end{array}$$
This yields
$$\frac{B_k(Z_{-k}(t))}{A_k(Z_{-k}(t))}\leq \frac{\pi_m\sigma_m^2}{1-Z_m(t)}\,.$$
As also $$ \pi_k\sigma_k^2\leq \frac{\pi_k\sigma_k^2}{1-Z_k(t)}\leq \frac{\pi_m\sigma_m^2}{1-Z_m(t)}\,,$$
we have thus proven that
$$\frac{\pi_k\sigma_k^2}{1-Z_k(t+1)}\leq M(t)\,.$$
This yields the monotonicity of $M(t)$ for $t\ge T$. Consequently,
$$\frac{1}{1-Z_j(t)}\leq \frac{M(T)}{\pi_j\sigma_j^2}\,,\qquad \forall j\in\mc V,\ t\ge T\,,$$
so that we have the a-priori bound
\be\label{apriori-bound}Z_j(t)\leq 1-\frac{\pi_j\sigma_j^2}{M(0)}\,,\qquad\forall j\in\mc V,\ t\ge T\,.\ee

We now consider the sequence $V(Z(t))$. By construction, it is non-increasing for $t\geq T$ and converges to some value $\bar V$ in $\R_+$. Consider now any limit point of the sequence $Z(t)$ and without loss of generality assume that $Z(t)\to \bar z$ (avoiding to formally pass to a subsequence). We know from previous considerations that $\bar z\in [0,1)^n$ and we want to show that $\bar z$ is a Nash equilibrium. Suppose not and let $i\in\mc \mc V$ be an agent that is not at equilibrium in configuration $\bar z$. Let $f_i: [0,1)^n\to  [0,1)^n$ be the function defined by
$f_i(z)=(\mc B_i(z_{-i}), z_{-i})$. We have $V(f_i(\bar z))=V(\bar z)-\epsilon$ for some $\epsilon >0$. Since $f_i$ is continuous, there exists $\delta>0$ such that $|z^1-z^2|<\delta$ yields
\be\label{contraction}\left|[V(f_i(z^1))-V(z^1)]- [V(f_i(z^2))-V( z^2)]\right|<\frac{\epsilon}{2}\ee
Consider now $\bar t\geq T$ such that $|Z(t)-\bar z|<\delta$ and $|V(Z(t))- V(\bar z)|<\epsilon/4$ for every $t\geq \bar t$ and pick any time $t\geq \bar t$ at which agent $i$ becomes active (with probability $1$ there are infinite such times).
From expression \eqref{contraction}, using the fact that $|Z(t)-\bar z|<\delta$ we deduce that
$V(f_i(Z(t)))-V(Z(t))<-\epsilon/2$.
Coupling now with $|V(Z(t))- V(\bar z)|<\epsilon/4$, we get that $V(f_i(Z(t)))-V(\bar z)<-\epsilon/4$ and this contradicts the fact that $V(Z(t))$ monotonically converges to $\bar z$. This proves that $\bar z$ is a Nash equilibrium. \qed
\end{IEEEproof}

We will see in examples below that the specific Nash equilibrium at which

\section{Numerical example}\label{subsec:exampl}

\begin{figure}[ht]
\centering
\includegraphics[width=0.5\columnwidth]{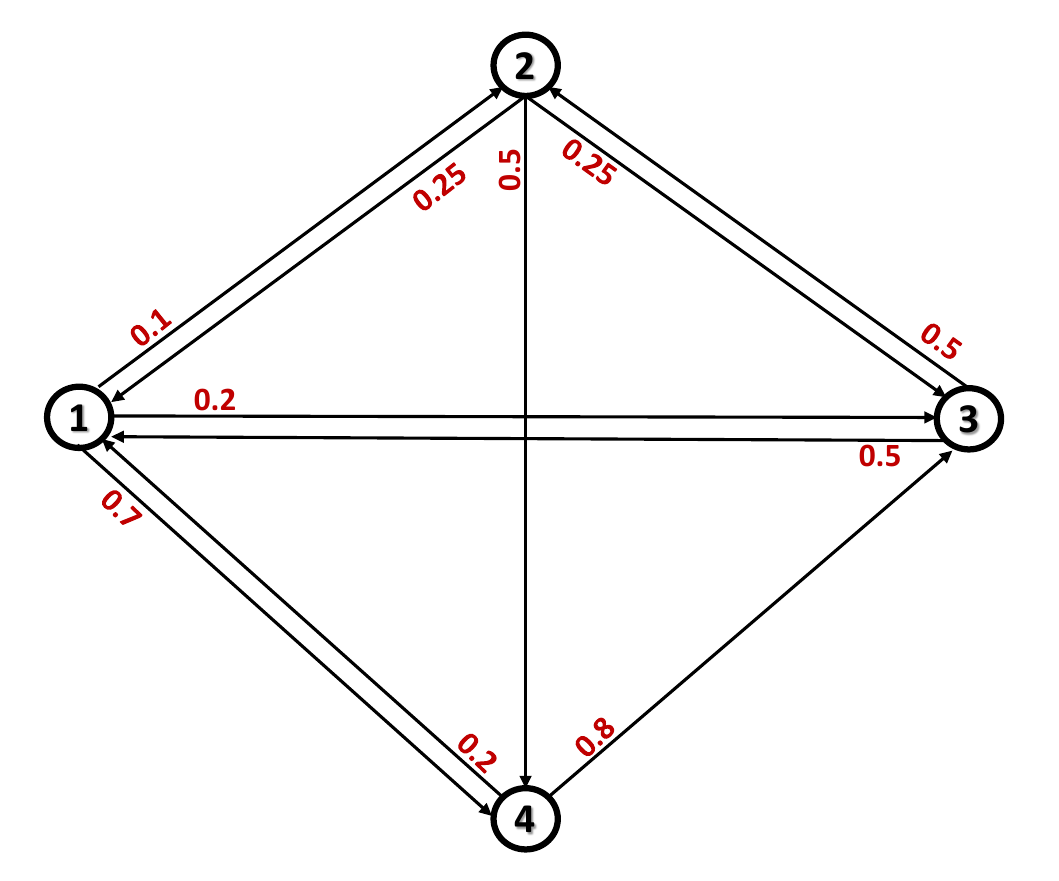}
\caption{Graph $\mc G_P$ of the matrix~\eqref{eq.p-example}}\label{fig.graph}
\end{figure}

Consider a social network with $n=4$ agents and influence matrix $P$ as follows
\begin{equation}\label{eq.p-example}
  P=\begin{bmatrix}
      0 & 0.1 & 0.2 & 0.7\\
      0.25 & 0 & 0.25 & 0.5\\
      0.5 & 0.5 & 0 & 0\\
      0.2 & 0 & 0.8 & 0
    \end{bmatrix}
\end{equation}
whose (weighted) graph is shown in Figure \ref{fig.graph}.
Let the vector of variances be
$$
\sigma^2=(0.32^2,\;0.35^2,\; 0.38^2,\; 0.29^2),
$$
i.e., agent $4$ is the ``wisest'' and agent $3$ is the most ``unwise''.
The centrality vector $\pi$ from~\eqref{pi} is found to be
$$\pi=(0.2507,\; 0.1783,\; 0.3064,\; 0.2646)\,.$$

If the agents had no self-influence ($z=0$), agent $3$ would be most influential, whereas agent $2$ would have the least social power among all agents.
For the Pareto-optimal profiles from the set $\mc P=\mc Z^*$, however, the highest self-confidence value belongs to agent $2$, who is neither the most ``central'' (in terms of eigenvector centrality) nor the wisest. In contrast, the most ``central'' agent $3$, has the smallest self-confidence value due to its high variance.
This follows from~\eqref{Z*}, since
$$[\pi]\sigma^2=(0.0257,\;0.0218,\;0.0442,\;0.0223).$$
In view of Corollary \ref{coro:no-nonstrct} all Nash equilibria are strict and coincide with the Pareto optimal self-confidence profiles.

We have simulated the asynchronous randomized best-response dynamics examined in Section~\ref{sec:dynam} for different initial conditions $z(0)$ and various pseudorandom activation signals\footnote{The pseudorandom sequences of active indices are generated using MATLAB's \texttt{randi} function, with the seed of the random number generator specified at the top of each plot.}. In all simulations, the self-confidence profile $Z(t)$ converges to one of the elements of the set $\mc Z^*$ (as defined in~\eqref{Z*}), i.e.,
\[
Z(t)\xrightarrow[t\to\infty]{} \one-\alpha[\pi]\sigma^2,
\]
where $\alpha\in(0,\alpha_*]$ (specified for each experiment) depends on $z(0)$ and the specific random sample path.

Our first experiment demonstrates that, even with the same initial condition, the final Nash equilibrium depends substantially on the sample path of the random activation signal.
This experiment also allows us to assess the conservatism of the a priori bound~\eqref{apriori-bound} derived in the proof of Theorem~\ref{theo:convergence}.
All agents start with self-confidence weights $z_i(0)=0.5$, which corresponds to $M(0)\approx 0.089$. The a priori bounds from~\eqref{apriori-bound} are, respectively,
$z_1\leq 0.71,\,z_2\leq 0.753,\,z_3\leq 0.5,\,z_4\leq 0.749$.

For the first random activation signal (corresponding to random generator seed 22, left subplot), $Z(t)$ converges to the vector $z^*\in\mc Z^*$, corresponding to $\alpha\approx 12.92$, or $z_1^*=0.67,z_2^*=0.72,z_3^*=0.43,z_4^*=0.71$, whose coordinates are thus sufficiently close to the maximal possible values for $z_i(t)$. Notice that all agents but for agent $3$ have increased their self-confidence weights.
The second random activation signal (random generator initialized with seed $1$) corresponds to a very different trajectory, resulting in the Nash equilibrium with $\alpha\approx 21.62$. Agents $2$ and $4$ slightly increase their self-confidence weights, and agent $1$'s weight is slightly decreased. Agent $3$'s self-confidence, however, decreases dramatically.
\begin{figure*}[htb]
	\centering
	\begin{subfigure}[b]{0.49\textwidth}
		\centering
		\includegraphics[width=\textwidth]{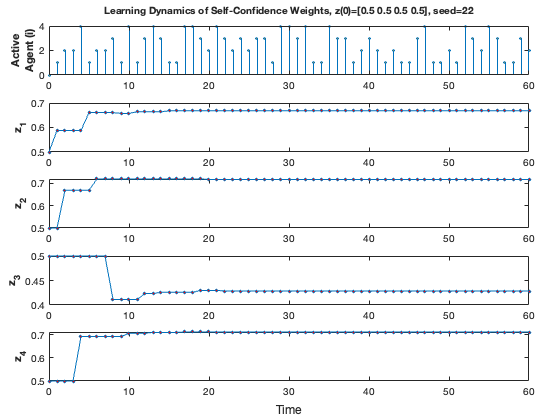}
		\caption{$\alpha\approx 12.92$}
		\label{fig:sub1}
	\end{subfigure}
	\hfill
	\begin{subfigure}[b]{0.49\textwidth}
		\centering
		\includegraphics[width=\textwidth]{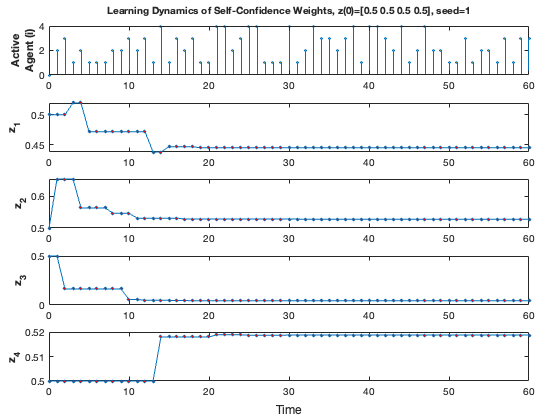}
		\caption{$\alpha\approx 21.62$}
		\label{fig:sub2}
	\end{subfigure}
	\caption{Convergence of the randomized learning dynamics for different activation signals}
	\label{fig:2}
\end{figure*}

Fig.~\ref{fig:2} illustrates the dependence between the initial self-confidence weight profile and the resulting Nash equilibrium,
using a fixed pseudo-random order of agent activation. In all experiments, agent $3$ (with the least wisdom) initially exhibits no self-confidence, whereas the other three agents are highly self-confident. Moreover, at least one agent is initially stubborn in every experiment.
In spite of similar initial conditions, the final Nash equilibria, convergence rates, and trajectories differ considerably. As expected, the parameter $\alpha$, which characterizes the final equilibrium, depends discontinuously on $z(0)$ when $\mc S(z(0)) \neq \emptyset$.
For instance, in the case where only agent $4$ is stubborn (top-right subplot), the final self-confidence weights of all agents exceed 0.9, and convergence to the equilibrium is very fast. In contrast, when only agent $2$ is stubborn, agent $3$'s final opinion is very small, and the evolution of $Z(t)$ is highly non-monotone—for example, agent $3$ acquires first a very large weight (close to yet strictly less than 1) and agent 4 even disconnects from the network by playing $z_4=1$, after which both self-weights drop rapidly.
\begin{figure*}[htb]
	\centering
	\begin{subfigure}[b]{0.49\textwidth}
		\centering
		\includegraphics[width=\textwidth]{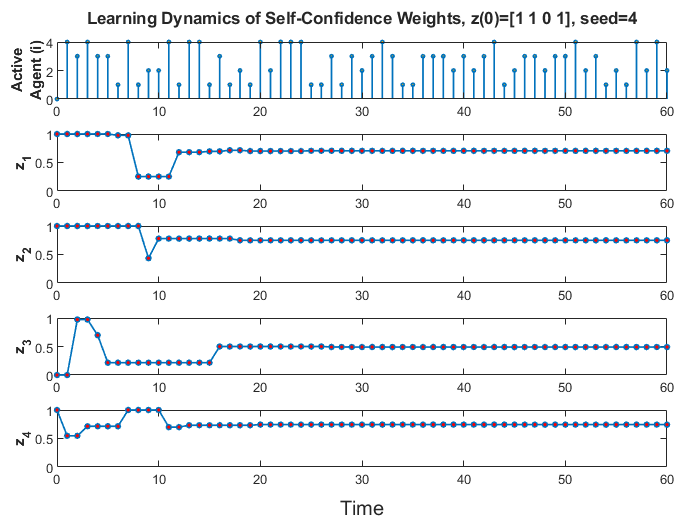}
		\caption{$\alpha\approx 11.51$}
		\label{fig:sub1}
	\end{subfigure}
	\hfill
	\begin{subfigure}[b]{0.49\textwidth}
		\centering
		\includegraphics[width=\textwidth]{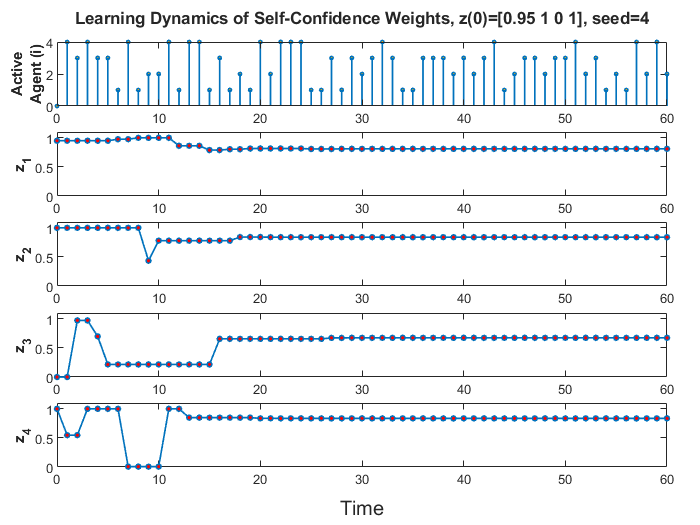}
		\caption{$\alpha\approx 7.39$}
		\label{fig:sub2}
	\end{subfigure}
	\begin{subfigure}[b]{0.49\textwidth}
		\centering
		\includegraphics[width=\textwidth]{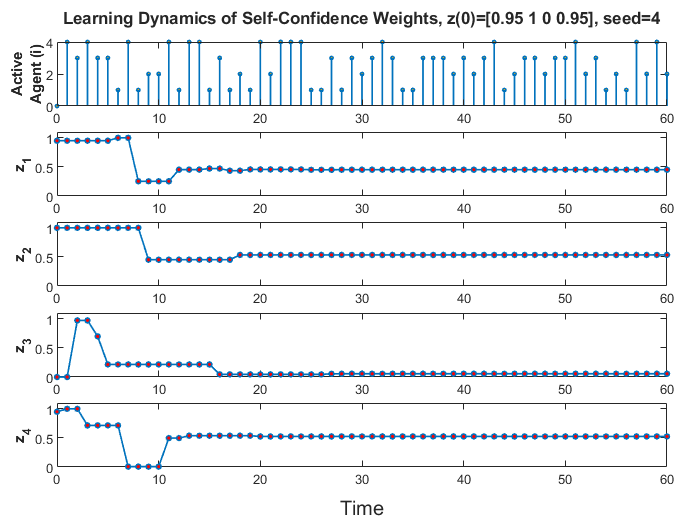}
		\caption{$\alpha\approx 21.33$}
		\label{fig:sub3}
	\end{subfigure}
	\begin{subfigure}[b]{0.49\textwidth}
		\centering
		\includegraphics[width=\textwidth]{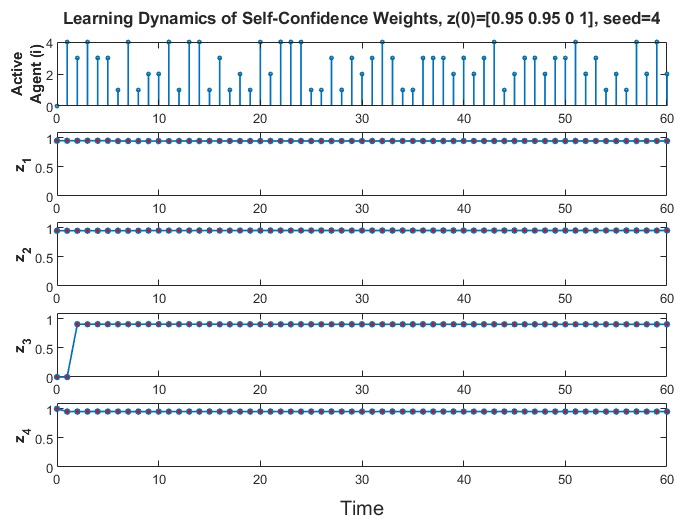}
		\caption{$\alpha\approx 2.12$}
		\label{fig:sub3}
	\end{subfigure}
	\caption{Convergence of the randomized learning dynamics for different $z(0)$}
	\label{fig:2}
\end{figure*}

\section{Conclusions and Future Works}
In this paper, we have considered a finite set of agents that first make uncorrelated unbiased noisy measurements of a common state of the world and then engage in a distributed averaging process, where every agent iteratively updates her estimate of the state of the world to a convex combination of her own current estimate and of a weighted average of those of her neighbors.  We have studied the network formation game resulting from letting each agent choose her self-confidence value ---corresponding to the weight assigned to her own estimate in the aforementioned convex combination--- with the aim of minimizing the variance of her asymptotic  estimate of the state of the world. For such network formation game, we have provided a full characterization of the set of pure strategy Nash equilibria and analyzed the asymptotic behavior of the asynchronous best response dynamics.

We believe that there are two main directions for further study worth being pursued. On the one hand, the problem can be generalized, e.g., by allowing the agents to choose not just their self-confidence value, but also to redistribute the relative influence weights assigned to their neighbors in the network, or by assuming that the agents' objective is not merely minimizing their asymptotic estimate's variance, but also, e.g., increasing their centrality in the resulting network~\citep{Castaldo.ea:2020,Catalano.ea:2024}. On the other hand, it would be worth investigating the behavior of other game-theoretic learning dynamics for the self-confidence evolution, beyond the asynchronous best-response dynamics.

\bibliographystyle{agsm}

\bibliography{mine,social}

\appendix
\section{Proof of Lemma \ref{lemma:limit}} \label{sec:proof-lemma:limit}

Let $\mc S=\mc S(z)$ denote the set of stubborn agents and $\mc R=\mc V\setminus\mc S$ be the complementary set of regular agents.
First, consider the case $\mc S=\emptyset$. In this case, since $0\le z<\one$ and $P$ is an irreducible aperiodic row-stochastic matrix, so is $W(z)=(I-[z])P+[z]$. It then follows from \cite[Theorem 1.8.3]{norris1998markov} that
$$H(z)=\lim_{t\to+\infty}W(z)^t=\one p(z)'\,,$$
where $ p(z)$ is the unique invariant probability distribution of $W(z)$. Hence, equation \ref{Hz} holds true.
Now, observe that
$$p(z)=W(z)'p(z)=((I-[z])P+[z])'p(z)\,,$$
if and only if
\be\label{ppi}P'[z]p(z)=[z]p(z)\,.\ee
Equation~\eqref{ppi} states that $[z]p(z)$ is an eigenvector of $P'$ associated to its leading eigenvalue $1$: since $P$ is irreducible, this implies that $[z]p(z)$ must be proportional to its unique invariant probability distribution $\pi$. This implies that Eq.~\eqref{def:p(z)} holds true with $\gamma(z)$ defined as in Eq.~\eqref{def:gamma(z)}. Hence, we have proven part (i) of the claim.

Now, consider the case $\mc S\ne\emptyset$.
Our proof in this case follows similar lines as the one of \cite[Theorem 2]{Como.Fagnani:2016}.
Decompose the stochastic matrix $W(z)$ as
\be\label{Wz-dec}W(z)=\left(\ba{cc} Q&B\\0&I\ea\right)\,,\ee
where $$Q=((I-[z])P+[z])_{\mc R\times\mc R}\,,$$ and
$$B=((I-[z])P)_{\mc R\times\mc S}\,.$$
Then, for every $t=0,1,\ldots$, we have
\be\label{Wzt}W(z)^t=\left(\ba{ccc} Q^t&\ \ \ds\sum_{0\le s<t}Q^sB\\0&I\ea\right)\,.\ee
Since $Q$ is a nonnegative square matrix, the Perron–Frobenius theorem
implies that there exits a nonnegative eigenvector $v\ne0$ such that
$$Q'v=\rho(Q)v\,,$$
where $\rho(Q)$ is the spectral radius of $Q$.
Let $\mc J=\{j\in\mc R:\,v_j>0\}$ be the support of $v$.
Since $\mc G_P$ is strongly connected, and $z_j<1$ for every $j$ in $\mc J$, there must exist some $j$ in $\mc J$ and $k$ in $\mc V\setminus\mc J$ such that
$W_{jk}(z)=(1-z_j)P_{jk}>0\,,$
so that
$$\sum_{i\in\mc J}Q_{ji}=\sum_{i\in\mc J}W_{ji}(z)\le1-W_{jk}(z)<1\,.$$
Hence,
$$
\ds\rho(Q)\sum_{i\in\mc J}v_i
=\ds\sum_{i\in\mc J}\sum_{j\in\mc R}Q_{ji}v_j
=\ds\sum_{j\in\mc J}v_j \sum_{i\in\mc J}Q_{ji}
<\ds\sum_{j\in\mc J}v_j\,,
$$
thus proving that $\rho(Q)<1$ (matrix $Q$ is Schur stable).
It follows that $I-Q$ is invertible and
\be\label{limQt} \lim_{t\to+\infty}Q^t=0\,,\qquad  \lim_{t\to+\infty}\sum_{0\le s<t}Q^s=(I-Q)^{-1}\,.\ee
It then follows from equations \eqref{Wzt} and \eqref{limQt} that
\be\label{Hz-dec}H(z)=\lim_{t\to+\infty}W(z)^t=\left(\ba{cc} 0\ \ &(I-Q)^{-1}B\\0&I\ea\right)\,.\ee
Notice that $H(z)$ is nonnegative since both matrices $(I-Q)^{-1}$ and $B$ are.
On the other hand, observe that $W(z)\one=\one$ implies that $Q\one+B\one=\one$, so that
$B\one=(I-Q)\one$, which in turn implies that
$$H(z)\one=\left(\!\!\ba{c}(I-Q)^{-1}B\one\\\one\ea\!\!\right)
=\left(\!\!\ba{c}(I-Q)^{-1}(I-Q)\one\\\one\ea\!\!\right)
=\one\,.$$
This proves that $H(z)$ is a row-stochastic matrix.

Now, observe that \eqref{Wz-dec} and \eqref{Hz-dec} imply that
\be\label{HRR}(W(z)H(z))_{\mc R\times\mc R}=0=(H(z))_{\mc R\times\mc R}\,,\ee
while
\be\label{HRS}\ba{rcl}(W(z)H(z))_{\mc R\times\mc S}&=& Q(I-Q)^{-1}B +B\\
&=&\left(Q(I-Q)^{-1}+I\right)B\\
&=&(I-Q)^{-1}B\\
&=&
(H(z))_{\mc R\times\mc S}\,,\ea\ee
where the second identity follows from the fact that
$$Q(I-Q)^{-1}+I=(I-Q)^{-1}\,,$$
since
$(Q(I-Q)^{-1}+I)(I-Q)=I-Q+Q=I\,.$
Equations \eqref{HRR} and \eqref{HRS} imply that, for every $i$ in $\mc R$ and $j$ in $\mc S$, we have
$$\ba{rcl}H_{ij}(z)
&=&\ds\sum_{k\in\mc V}W_{ik}(z)H_{kj}(z)\\
&=&\ds(1-z_i)\sum_{k\in\mc V}P_{ik}H_{kj}(z)+z_iH_{ij}(z)\,.
\ea$$
Since $1-z_i>0$ for every $i$ in $\mc R$, the above implies that the first equation in \eqref{lin-sys} holds true. That the second equation in \eqref{lin-sys} holds true is immediate from \eqref{Hz-dec}. Finally, uniqueness of the solution of the linear system \eqref{lin-sys} follows from invertibility of $I-Q$, thus completing the proof of part (ii) of the claim.
\qed

\section{Proof of Lemma \ref{lemma:min}}\label{sec:proof-lemma:min}
First, notice that, when $\one'\mu=1$, one has
\[
\sum_{i\in\mc V}(\mu_i\sigma_i)\sigma_i^{-1}=\sum_{i\in\mc V}\mu_i=1.
\]
By applying the Cauchy-Schwarz inequality to the identity above, one gets that
\[
\sum_{i\in\mc V}(\mu_i\sigma_i)^2\sum_{i\in\mc V}\sigma_i^{-2}\geq 1,
\]
where the inequality is strict unless the two vectors are collinear, i.e., $\mu_i=\alpha\sigma_i^{-2}$ for some $\alpha$ in $\mathbb{R}$.  The latter condition is equivalent to~\eqref{mu}, recalling that $\one'\mu=1$.\qed

\end{document}